\documentclass[11pt]{amsart}

\usepackage{amssymb,latexsym}

\usepackage{enumerate}

\makeatletter

\@namedef{subjclassname@2010}{%

  \textup{2010} Mathematics Subject Classification}

\makeatother

\numberwithin{equation}{section}

\begin{document}

%\baselineskip=17pt 
%arXiv version has this option "on"

\title[Focal Rigidity of Flat Tori]{Focal Rigidity of Flat Tori}

\author[F. H. Kwakkel]{Ferry Kwakkel}

\address{Instituto de Matematica Pura e Aplicada \\ Rio de Janeiro \\ Brazil}

\email{kwakkel@impa.br}

\author[M. Martens]{Marco Martens}

\address{SUNY at Stony Brook \\ New York \\ United States}

\email{marco@math.sunysb.edu}

\author[M. M. Peixoto]{Mauricio Peixoto}

\address{Instituto de Matematica Pura e Aplicada \\ Rio de Janeiro \\ Brazil}

\email{peixoto@impa.br}

\date{}

\newtheorem{thm}{Theorem}
\newtheorem{cor}[thm]{Corollary}
\newtheorem{lem}[thm]{Lemma}
\newtheorem{prop}[thm]{Proposition}
\newtheorem{notation}[thm]{Notation}

\newtheorem*{thma}{Theorem A}
\newtheorem*{strat}{Existence Theorem}
\newtheorem*{mostow}{Mostow's Rigidity Theorem}
\newtheorem*{open}{Open problem}

\theoremstyle{definition}

\newtheorem{defn}{Definition}
\newtheorem{rem}{Remark}

\theoremstyle{open}

\newtheorem{quest}{Question}

\numberwithin{equation}{section}
\newcommand{\figref}[1]{Figure~\ref{#1}}
\newcommand{\pichere}[2]
{\begin{center}\includegraphics[width=#1\textwidth]{#2}\end{center}}

\newcommand{\bignote}[1]{\begin{quote} \sf #1 \end{quote}}

\newcommand{\QED}{\rlap{$\sqcup$}$\sqcap$\smallskip}

\renewcommand{\Im}{\operatorname{Im}}
\renewcommand{\exp}{\operatorname{exp}}

\def\sss{\subsubsection}

\newcommand{\di}{\partial}
\newcommand{\dibar}{\bar\partial}
\newcommand{\hookra}{\hookrightarrow}
\newcommand{\ra}{\rightarrow}
\newcommand{\hra}{\hookrightarrow}
\newcommand{\imply}{\Rightarrow}
\def\lra{\longrightarrow}
\newcommand{\wc}{\underset{w}{\to}}
\newcommand{\tu}{\textup}
\newcommand{\SO}{\textup{SO}}
\newcommand{\Int}{\textup{Int}}

\newcommand{\eps}{{\varepsilon}}
\newcommand{\epsi}{{\epsilon}}
\newcommand{\veps}{{\varepsilon}}
\newcommand{\De}{{\Delta}}
\newcommand{\de}{{\delta}}
\newcommand{\la}{{\lambda}}
\newcommand{\La}{{\Lambda}}
\newcommand{\si}{{\sigma}}
\newcommand{\Si}{{\Sigma}}
\newcommand{\Om}{{\Omega}}
\newcommand{\om}{{\omega}}

\newcommand{\AAA}{{\mathcal A}}
\newcommand{\BB}{{\mathcal B}}
\newcommand{\CC}{{\mathcal C}}
\newcommand{\DD}{{\mathcal D}}
\newcommand{\EE}{{\mathcal E}}
\newcommand{\EEE}{{\mathcal O}}
\newcommand{\II}{{\mathcal I}}
\newcommand{\FF}{{\mathcal F}}
\newcommand{\GG}{{\mathcal G}}
\newcommand{\JJ}{{\mathcal J}}
\newcommand{\HH}{{\mathcal H}}
\newcommand{\KK}{{\mathcal K}}
\newcommand{\LL}{{\mathcal L}}
\newcommand{\MM}{{\mathcal M}}
\newcommand{\NN}{{\mathcal N}}
\newcommand{\OO}{{\mathcal O}}
\newcommand{\PP}{{\mathcal P}}
\newcommand{\QQ}{{\mathcal Q}}
\newcommand{\QM}{{\mathcal QM}}
\newcommand{\QP}{{\mathcal QP}}
\newcommand{\QL}{{\mathcal Q}}

\newcommand{\RR}{{\mathcal R}}
\renewcommand{\SS}{{\mathcal S}}
\newcommand{\TT}{{\mathcal T}}
\newcommand{\TTT}{{\mathcal P}}
\newcommand{\UU}{{\mathcal U}}
\newcommand{\VV}{{\mathcal V}}
\newcommand{\WW}{{\mathcal W}}
\newcommand{\XX}{{\mathcal X}}
\newcommand{\YY}{{\mathcal Y}}
\newcommand{\ZZ}{{\mathcal Z}}

\newcommand{\A}{{\mathbb A}}
\newcommand{\C}{{\mathbb C}}
\newcommand{\bC}{{\bar{\mathbb C}}}
\newcommand{\D}{{\mathbb D}}
\newcommand{\Hyp}{{\mathbb H}}
\newcommand{\J}{{\mathbb J}}
\newcommand{\Ll}{{\mathbb L}}
\renewcommand{\L}{{\mathbb L}}
\newcommand{\M}{{\mathbb M}}
\newcommand{\N}{{\mathbb N}}
\newcommand{\Q}{{\mathbb Q}}
\newcommand{\R}{{\mathbb R}}
\newcommand{\B}{{\mathbb B}}
\newcommand{\T}{{\mathbb T}}
\newcommand{\V}{{\mathbb V}}
\newcommand{\U}{{\mathbb U}}
\newcommand{\W}{{\mathbb W}}
\newcommand{\X}{{\mathbb X}}
\newcommand{\Z}{{\mathbb Z}}
\newcommand{\ball}{{\mathbb B}}
\newcommand{\Id}{\textup{Id}}
\newcommand{\cl}{\textup{Cl}}
\newcommand{\GL}{\textup{GL}}

\newcommand{\VVV}{{\mathbf U}}
\newcommand{\UUU}{{\mathbf U}}

\newcommand{\tT}{{\mathrm{T}}}
\newcommand{\tD}{{D}}
\newcommand{\hyp}{{\mathrm{hyp}}}

\newcommand{\f}{{\bf f}}
\newcommand{\g}{{\bf g}}
\newcommand{\h}{{\bf h}}
\renewcommand{\i}{{\bar i}}
\renewcommand{\j}{{\bar j}}

\catcode`\@=12

\def\Empty{}
\newcommand\oplabel[1]{
  \def\OpArg{#1} \ifx \OpArg\Empty {} \else
   \label{#1}
  \fi}

%%%%%%%%%%%%%%%%%%%%%%%%%%%%%%%%%%%%%%%%%%%%%%%%%%%%%%%%%%%%%%%%%%%%%
% Insert a postscript figure using psfig.
% Usage: \realfig{label}{filename}{caption}
%
% uses psfig macros: must have \input{psfig} in the preamble to use
% it. 
%%%%%%%%%%%%%%%%%%%%%%%%%%%%%%%%%%%%%%%%%%%%%%%%%%%%%%%%%%%%%%%%%%%%%

\long\def\realfig#1#2#3#4{
\begin{figure}[htbp]
%%%\centerline{\psfig{figure=#3,height=#2}}
\centerline{\psfig{figure=#2,width=#4}}
\caption[#1]{#3}
\oplabel{#1}
\end{figure}}

\numberwithin{figure}{section}

%&&&&&&&&&&&& List of figures &&&&&&&&&
%
%&&&&&&&&&&&&&&&&&&&&&&&&&&&&&&&&&&&&&&&&&&&&&&&&&&&&&

\newcommand{\comm}[1]{}
\newcommand{\comment}[1]{}

\begin{abstract}
Given a closed Riemannian manifold $(M, g)$, there is a partition of its tangent bundle $TM = \bigcup_i \Sigma_i$, called the focal decomposition of $TM$. The sets $\Sigma_i$ are closely associated to focusing of geodesics of $(M,g)$, i.e. to the situation 
where there are exactly $i$ geodesic arcs of the same length joining points $p$ and $q$ in $M$. In this note, we study the topological structure of the focal decomposition of a closed Riemannian manifold and its relation with the metric structure of the manifold. 
Our main result is that the flat $n$-tori are focally rigid, in the sense that if two flat tori are focally equivalent, then the tori are isometric up to rescaling. 
\end{abstract}



\maketitle

\section{Definitions and Statement of Results}\label{sec_intro}

In general, topological characteristics of a Riemannian manifold do not determine its geometry, that is, its metric structure. 
However, starting in the 1960s, examples have been discovered for which such characteristics do determine the geometry. 
The manifolds can not be deformed without changing the characteristic. One speaks of rigidity. The prototype rigidity result is due to Mostow~\cite{Mos}.

\begin{mostow}\label{thm_mostow}
Two compact hyperbolic $n$-manifolds, with $n \geq 3$, with isomorphic fundamental groups are isometric.
\end{mostow}

Given an analytic manifold $M$, we denote $\RR^{\omega}(M)$ the class of analytic Riemannian structures on $M$. This class of structures is dense in the class $\RR^{\infty}(M)$ of smooth Riemannian structures in the $C^{\infty}$ strong Whitney topology, see~\cite{hirsch}.
Let $(M, g)$ and $(\widetilde{M}, \widetilde{g})$ be closed (compact and boundaryless) analytic manifolds. Two Riemannian manifolds are said to be {\em isometric up to a rescaling} if, up to a constant rescaling of the metric $\widetilde{g}$, 
the manifolds $(M, g)$ and $(\widetilde{M}, \widetilde{g})$ are isometric. 

In~\cite{pugh}, the notion of {\em focal stability} is introduced for smooth Riemannian manifolds in dimension $n$. In analogy with structural stability for dynamical systems, this gives rise to the focal stability conjecture: given $p \in M$, 
the generic Riemannian structure is focally stable at $p$. In~\cite{pugh_2}, if $n=2$ and there are no conjugate points, the above conjecture has been shown to be true. The notion of focal stability in~\cite{pugh} is {\em local}.
The purpose of this note is to give the first examples of closed manifolds, namely flat $n$-tori, and show that the focal decomposition is {\em globally} rigid. The concept of focal rigidity and related concepts were first studied for two-tori in~\cite{kwa}.

\begin{thma}
Focally equivalent flat $n$-tori, where $n \geq 2$, are isometric, up to rescaling. 
\end{thma}

The proof, though elementary, makes use of the notion of Brillouin zones intrinsic to the focal decomposition, and will also be shown in a forthcoming note~\cite{ferry}, in which the focal rigidity of hyperbolic manifolds is considered. 
In order to define the notion of {\em focal equivalence}, we first need some preliminary definitions. Fix a complete $C^{\omega}$ manifold $(M, g)$.

\begin{defn}\label{defn_index}
The focal index, $I(p,v)$, of the vector $v \in T_pM$ is defined by 
\begin{equation*}
I(p,v) = \# \left\{ w \in T_p M ~|~ |v| = |w| ~ \tu{and} \exp_p(v) = \exp_p(w) \right\}.
\end{equation*}
The {\it focal component} of index $i$ at $p$ is
\begin{equation*}
\sigma_i (p) = \left\{ v \in T_p M ~|~ i = I(p,v) \right\}.
\end{equation*}
\end{defn}

Vectors $v \in \sigma_i(p)$ are equivalent modulo exponentiation to $i-1$ other vectors of $T_p M$ of equal length. 

\begin{defn}[Focal decomposition]\label{defn_decomp}
The partition of $T_p M$ into its focal components $\{ \sigma_i(p) \}_{i=1}^{\infty}$ is called its focal decomposition at $p$; we have
\begin{equation}\label{eq_decomp}
T_p M = \bigcup_{i=1}^{\infty} \sigma_i(p) ~\textup{and}~ \sigma_i(p) \cap \sigma_j(p) = \emptyset, ~\textup{if}~ i \neq j.
\end{equation}
The tangent bundle has a corresponding focal decomposition $\{ \Sigma_i\}_{ i=1}^{\infty}$, where $\Sigma_i = \bigcup_{p \in M} \sigma_i (p)$ 
and correspondingly
\begin{equation}\label{eq_decomp_2}
TM = \bigcup_{i=1}^{\infty} \Sigma_i ~\textup{and}~ \Sigma_i \cap \Sigma_j = \emptyset, ~\textup{if}~ i \neq j.  
\end{equation}
\end{defn}

The focal decomposition only depends on the Riemannian metric $g$. It is also a global concept: all geodesics passing through $p$ play a role in the construction of $\sigma_i$ and all geodesics in the manifold $M$ play a role in the construction of 
the sets $\Sigma_i$. We will abbreviate $\{ \sigma_i \}_{i=1}^{\infty} = \{ \sigma_i\}$ and $\{ \Sigma_i \}_{i=1}^{\infty} = \{ \Sigma_i\}$.

\begin{defn}[Focal equivalence]\label{focal_equiv}
Two complete analytic manifolds $(M_1, g_1)$ and $(M_2, g_2)$ are {\em focally equivalent}, if there exists an orientation-preserving homeomorphism $\varphi \colon TM_1 \ra TM_2$, with $\psi \colon M_1 \ra M_2$ the homeomorphism induced on the zero section,
such that for every $p \in M_1$ and $q = \psi(p) \in M_2$,
\begin{enumerate}
\item[\tu{(i)}] $\varphi_{|_{T_pM_1}} \colon T_pM_1 \ra T_{q} M_2$ and $\varphi_{|_{T_pM_1}}(0) = 0$,
\item[\tu{(ii)}] $\varphi_{|_{T_pM_1}} (\sigma^1_i(p)) = \sigma^2_i(q)$, for every $1 \leq i \leq \infty$.
\end{enumerate}
\end{defn}

It follows from (i) and (ii) in Definition~\ref{focal_equiv} that $\varphi(\Sigma_i^1) = \Sigma_i^2$. It is verified that focal equivalence indeed defines an equivalence relation. 
Further, manifolds that are isometric, up to rescaling, are focally equivalent. 

\section{The Focal Decomposition and Brillouin Zones}\label{sec_foc_decomp}

The concept of the focal decomposition was introduced in~\cite{kupka} under the name of {\it sigma decomposition} 
in the context of the two-point boundary value problem for ordinary second order differential equations. 
It was then naturally extended to complete Riemannian manifolds of any dimension by I. Kupka and M. Peixoto in~\cite{KP} 
as outlined in the introduction.

In~\cite{pugh}, a survey is given of the focal decomposition of complete smooth ($C^{\infty}$) Riemannian manifolds. The main results are that for a generic\footnote{containing  a residual subset in the space of smooth metrics $\mathcal{R}^{\infty}(M)$ equipped with 
the $C^{\infty}$ strong Whitney topology.} Riemannian metric, the index $i$ of $\sigma_i$ is bounded by $\dim(M) + 1$ and the index $i$ of $\Sigma_i$ is bounded by $2 \dim(M) + 2$. Also in~\cite{pugh}, it is indicated how the concept of focal decomposition 
relates naturally to physical concepts such as the Brillouin zones of a crystal, the semi-classical quantization via Feynman path integrals and to Diophantine equations in number theory.

Recall that a Riemannian manifold $(M,g)$ is said to be complete if for every point $p \in M$ the exponential map at $p$, $\exp_p : T_p M \rightarrow M$, is defined at all points of $T_p M$. Every compact Riemannian manifold is complete by the 
Hopf-Rinow Theorem. Even though the focal decomposition of a smooth Riemannian manifold always exists, the choice of putting the discussion in the setting of analytic Riemannian manifolds is motivated by the following, see~\cite{KP}.

\begin{strat}\label{thm_strat} 
Let $(M,g)$ be a complete analytic Riemannian manifold, then there is an analytic Whitney stratification of $TM$ such that each $\Sigma_i$ is the union of strata of this stratification. Similarly for the restricted problem with base point $p \in M$: there 
is an analytic Whitney stratification of $T_pM$ such that every set $\sigma_i$ is the union of strata of this stratification. Analytically stratifiable means that the focal components $\sigma_i$ can be expressed as locally finite disjoint unions of strata. 
\end{strat}

\begin{rem}\label{rem_strat}
This result is sharp in the sense that if the Riemannian metric is smooth but not analytic, the focal components may be non-stratifiable, i.e. topologically pathological, see~\cite[page 39]{pugh}.
\end{rem}

\begin{rem}
It follows from the Angle Lemma in~\cite{KP} that no $\sigma_i$ with $i \geq 2$ has any interior points. 
Combined with the Whitney stratification property, one can show that $\sigma_1$ has full measure.
\end{rem}

Closely related to the focal decomposition, is the notion of {\it Brillouin zones}~\cite{pei2}, which in a way organizes the structure of the focal decomposition at a point $p$.

\begin{defn}[Brillouin Zones]\label{defn_bril}
For $v \in T_p M$, we define the Brillouin index
\begin{equation}\label{eq_bril}
B(p,v) =  \# \left\{ w \in T_p M ~|~ |w | \leq |v|, ~ \exp_p(w) = \exp_p(v) \right\}.
\end{equation}
For every integer $k \geq 1$, the $k$-th Brillouin zone is the interior $\Int(B_k(p))$, of the set $B_k(p) = \{ v \in T_p M ~|~ B(p,v) = k \}$ of all points with Brillouin index $k$.
\end{defn}

\begin{rem}
The $k$-th Brillouin zone is also defined to be $B_k(p)$, rather than its interior. The significance in taking the interior in our definition is that, in general, $\sigma_1$ need not be open, cf.~\cite[page 37]{pugh}.
\end{rem}

The first Brillouin zone appears in many different places; for example, it appears as the \emph{Wigner-Seitz cell} in physics, as the \emph{Dirichlet region} of a Fuchsian group in geometry and as the \emph{Voronoi cell} in the study of 
circle packings, see~\cite{jones} and~\cite{skri}. 

\section{Focal Rigidity of Flat Tori}\label{sec_foc_rig}

As the flat torus admits a transitive group of isometries, the focal decomposition is identical at every basepoint. Therefore, it suffices to study the focal decomposition at a single base point. 
Given a flat torus $(\T^n, g)$, it is isometric to a torus of the form $(\R^n \slash \Gamma, g_{\tu{can}})$, where $\Gamma$ is a translation group of rank $n$ and $g_{\tu{can}}$ the canonical metric induced from the Euclidean metric
of $\R^n$, see~\cite{wolf}. Therefore, we may as well assume that 
\begin{equation}
\T^n = \R^n \slash \Gamma \quad \tu{and} \quad \widetilde{\T}^n = \R^n \slash \widetilde{\Gamma}, 
\end{equation}
both equipped with the induced Euclidean metric $g_{\tu{can}}$, where $\Gamma$ and $\widetilde{\Gamma}$ are translation groups of rank $n$ acting on $\R^n$. In what follows, we identify the tangent 
plane $T_p\T^n$ with the universal cover $\R^n$ of $\T^n$ and we denote $d$ the distance function on $\R^n$ induced by the standard Euclidean metric. Choosing the basepoints to be $0 \in \T^n$ and $0 \in \widetilde{\T}^n$, where $0 = \pi(0)$ with $0 \in \R^n$ the 
origin, we denote $\sigma_i(0) = \sigma_i$ and $B_k(0) = B_k$ for brevity (and similarly for $\widetilde{\T}^n$). In terms of this uniformization and notation, the focal decomposition and corresponding Brillouin zones have the following description. 

\begin{defn}\label{bril_plane}
Given a flat torus $\R^n \slash \Gamma$, let $\Lambda = \Gamma(0)$ the associated lattice. The \emph{Brillouin hyperplane} $V_{\lambda} \subset \R^n$, $\lambda \in \Lambda \subset \R^n$, is defined to be the set 
$V_{\lambda} = \{v \in \R^n ~ \vert ~ d(0,v) = d(v, \lambda) \}$. For $0 \in \Lambda$, we define $V_0 = \emptyset$ and we denote $\mathcal{V} = \bigcup_{{\lambda} \in \Lambda} V_{\lambda}$.
\end{defn}

In what follows, we will refer to a Brillouin hyperplane simply as B-plane.

\begin{rem}\label{rem_loc_finite}
As the lattice $\Lambda$ is discrete, the set of B-planes $\mathcal{V} \subset \R^n$ is {\it locally finite} in the sense that for every 
$v \in \mathcal{V}$ and every compact set $K \subset \R^n$ containing $v$, only finitely many B-planes meet $K$.
\end{rem}

Let $\ell_v \subset \R^n$ be the open (i.e. not containing $0$ and $v$) line segment connecting $0$ and $v$ in $\R^n$ and define the indices
\[ \iota(v) = \# \left\{ \lambda \in \Lambda ~\vert~ V_{\lambda} \cap \ell_v \neq \emptyset \right\}  \quad \textup{and} \quad 
\mu(v) = \# \left\{ {\lambda} \in \Lambda ~ \vert ~ V_{\lambda} \ni v \right\}. \]

The following result gives a description of the focal decomposition and Brillouin zones for the torus, 
see~\cite{pei2} and~\cite{veerman}.

\begin{prop}\label{prop_focal_torus}
The focal decomposition of the torus $\T^n$ (at the basepoint $0$), is given by the following.
\begin{itemize}
\item[\tu{(1)}] $\sigma_i = \left\{ v \in \R^n ~|~ \mu(v) = i - 1 \right\}$ and, consequently, $\R^n \setminus \sigma_1 = \mathcal{V}$. 
\item[\tu{(2)}] For $v \in \R^n$, $v \in \Int(B_k)$ if and only if $\iota(v) = k -1 $ and $\mu(v) = 0$. 
\item[\tu{(3)}] $\cl( \Int (B_k)) = B_k$.
\end{itemize}
\end{prop}

Let us now turn to the proof of Theorem A. We denote by $\varphi : \R^n \rightarrow \R^n$ the orientation preserving homeomorphism for which $\varphi(0) = 0$ and $\varphi(\sigma_i) = \widetilde{\sigma}_i$,
which is the restriction to a single tangent plane of the homeomorphism $\varphi$ defined globally between the tangent bundles of the tori.
Denote $\mathcal{V} = \bigcup_{\lambda \in \Lambda} V_{\lambda}$ and $\widetilde{\mathcal{V}} = \bigcup_{\tilde{\lambda} \in \widetilde{\Lambda}} V_{\tilde{\lambda}}$ the union of all B-planes. 
Let $v \in \R^n$ and let $\{ V_{\lambda_j} \}_{j=1}^m$ be the set of B-planes for which $v \in V_{\lambda_j}$, $1 \leq j \leq m$. Note that $m= \mu(v)$. Let $\nu(v)$ be the dimension 
of the span of the vectors $\lambda_1,...,\lambda_m$. Further, a B-plane $V_{\lambda}$ separates $\R^n$ into two connected components, i.e.
\begin{equation}
\R^n - V_{\lambda} = \HH_1 \cup \HH_2, \quad \tu{with} ~ \HH_1 \cap \HH_2 = \emptyset.
\end{equation}
We say that $V_{\lambda}$ \emph{separates} the points $v$ and $w$, if $v \in \HH_1$ and $w \in \HH_2$. 

\begin{lem}[Plane Lemma]\label{planes}
For every $\lambda \in \Lambda$, there exists a unique $\widetilde{\lambda} \in \widetilde{\Lambda}$ such that $\varphi(V_{\lambda}) = V_{\tilde{{\lambda}}}$.
\end{lem}

In order to prove this lemma, we use the following auxiliary lemmas. 

\begin{lem}\label{sublemma_1}
Let $V = \R^k$ be a Euclidean space and $W \subset V$ a closed codimension two subset. 
Then $V \setminus W$ is path-connected and $V \setminus W$ is dense in $V$.
\end{lem}

\begin{proof}
As $W \subset V$ is closed, $V \setminus W$ is open. Further,  as $W$ has codimension two relative to $V$, it is clear that $V \setminus W$ is dense in $V$.
We need to show that $V \setminus W$ is path-connected. To this end, take any two points $x,y \in V \setminus W$ and consider the straight segment $\ell \subset V$ connecting $x$ and $y$. 
Take a product neighborhood $N$ of the segment $\ell$. The set of segments in $N$ that meet $W$ has codimension one in $N$. Therefore, we can find a segment $\ell'$ arbitrarily close to $\ell$ such that $\ell' \cap W = \emptyset$.
Further, since $x,y \in V \setminus W$ and $V\setminus W$ is open, by choosing the segment $\ell'$ close enough to $\ell$, we can connect $x$ with $x'$ and $y$ with $y'$ by an arc in $V\setminus W$, thus producing the desired path from $x$ to $y$. 
\end{proof}

\begin{lem}\label{sublemma_2}
Let $V_1, ... , V_m \subset \R^n$ and $V_1',...,V_m' \subset \R^n$ be two collections of codimension one hyperplanes in $\R^n$, where $m \geq 2$, such that 
\begin{enumerate}
\item[\tu{(1)}] $V_k$ and $V_k'$ pass through $0 \in \R^n$, for all $1 \leq k \leq m$, and
\item[\tu{(2)}] the dimension of the span of the normal vectors to $V_1, ... , V_m$ and $V_1',...,V_m'$ is both two.
\end{enumerate}
Let $h \colon \R^n \ra \R^n$ be a homeomorphism such that $h(0) = 0$ and $h(\bigcup_{k=1}^m V_k) = \bigcup_{k=1}^m V_k'$.
Then, after a suitable relabeling if necessary, we have $h(V_k) = V_k'$.
\end{lem}

\begin{proof}
Since the dimension of the span of the two collections of hyperplanes passing through the origin is two, after pre- and postcomposing with suitable rotations, we may assume that 
\begin{equation}
\bigcup_{k=1}^m V_k = Q \times \R^{n-2}  \quad \tu{and} \quad \bigcup_{k=1}^m V_k' = Q' \times \R^{n-2}, 
\end{equation}
where $Q \subset \R^2$ is a union of $m$ straight Euclidean lines $L_k$ passing through the origin $0 \in \R^2$, according to the rule $V_k = L_k \times \R^{n-2}$, where each line $L_k$ is perpendicular in $\R^2$ to the normal vector to the plane $V_k$, where $k=1, ..., m$, 
one line $L_k$ precisely corresponding  to a plane $V_k$. Similarly for $Q'$. Thus $\R^n \setminus \bigcup_{k=1}^m V_k$ consists of $2m$ connected components and so does $\R^n \setminus \bigcup_{k=1}^m V_k'$. Take a simple closed curve $\gamma \subset \R^2 \times \{ 0 \}$ winding around the origin $0 \in \R^2$, such that the curve $\gamma$ intersects every line $L_k$ exactly once. Label the intersection points $\{ z_s \}_{s=1}^{2m}$ of $\gamma$ with $\bigcup_{k=1}^m L_k \setminus \{ 0 \}$ according to their cyclic ordering along $\gamma$, relative to an orientation and initial point $z_1$, which we may assume to be contained in $L_1$. Further, we may assume that $h(z_1) \in V_1'$. In this labeling, $z_s$ and $z_{s+m}$ belong to the same line $L_k$ and thus to the same plane $V_k$. Given a line $L_k$, the B-plane $V_k$ corresponding to $L_k$ separates $\R^n$
into two connected components $\HH_k^+$ and $\HH_k^-$, each containing precisely $m$ connected components of the complement of $\R^n \setminus \bigcup_{k=1}^m V_k$. Since $h$ sends connected components of $\R^n \setminus \bigcup_{k=1}^m V_k$, in a one-to-one correspondence,
to the connected components of $\R^n \setminus \bigcup_{k=1}^m V_k'$, the images $h(\HH_k^+)$ and $h(\HH_k^-)$ contain precisely $m$ connected components of $\R^n \setminus \bigcup_{k=1}^m V_k'$ each. This can only be if $h(z_s)$ and $h(z_{s+m})$ belong to the same image 
plane. It thus follows, upon an appropriate labeling, that $h(V_k) = V_{k}'$, as required.
\end{proof}

\begin{proof}[Proof of Lemma~\ref{planes}]
Fix a B-plane $V_{\lambda}$. Let $v \in V_{\lambda}$ for which $\mu(v) = 1$ and let $\varphi(v) = \tilde{v}$ and therefore $\mu(\tilde{v}) = 1$. Let $V_{\tilde{{\lambda}}}$ be the unique B-plane passing through $\tilde{v}$. It suffices to show 
that $\varphi(V_{\lambda}) \subseteq V_{\tilde{{\lambda}}}$. Indeed, if this is proved, then by symmetry $\varphi^{-1}(V_{\tilde{{\lambda}}}) \subseteq V_{\lambda}$ and therefore $\varphi(V_{\lambda}) = V_{\tilde{{\lambda}}}$. 
To prove that $\varphi(V_{\lambda}) \subseteq V_{\tilde{{\lambda}}}$, consider 
\[ \mathcal{I}_{\lambda} := \{ v \in V_{\lambda} ~|~ \mu(v) \geq 2, \nu(v) \geq 3\} \subset V_{\lambda}, \]
and define 
\begin{equation}
W_{\lambda} = \mathcal{I}_{\lambda} \cup \left( \bigcup_{\widetilde{\lambda} \in \widetilde{\Lambda}} \varphi^{-1}(\mathcal{I}_{\widetilde{\lambda}}) \cap V_{\lambda} \right) \subset V_{\lambda}
\end{equation}
As $\mathcal{I}_{\lambda}$ and $\varphi^{-1}(\mathcal{I}_{\widetilde{\lambda}}) \cap V_{\lambda}$ have codimension two relative to $V_{\lambda}$, are closed subsets, and since the union of B-planes relative to $\widetilde{\Lambda}$ is a locally finite union of planes, 
by Lemma~\ref{sublemma_1}, it follows that $V_{\lambda} \setminus W_{\lambda}$ is open and path-connected and $V_{\lambda} \setminus W_{\lambda}$ is dense in $V_{\lambda}$, since a locally finite union of closed codimension two subsets 
is closed and has codimension two. Furthermore, $V_{\lambda} \setminus W_{\lambda}$ has the property that every point $v \subset V_{\lambda} \setminus W_{\lambda}$, for which $\mu(v) = m \geq 2$, we have $\nu(v) = 2$ and furthermore, $\mu(\tilde{v}) =m$ and $\nu(\tilde{v}) = 2$. 
Since $V_{\lambda} \setminus W_{\lambda}$ is open, we can take a small Euclidean ball $U_v$ centered at $v$, and by shrinking the ball if necessary, we may assume that $U_v$ only meets those B-planes that 
pass through $v$. Since $\varphi$ is a homeomorphism, $\varphi(U_v)$ will only meet those B-planes that pass through $\tilde{v}$. Moreover, as $U_v$ is an open ball, which is homeomorphic to $\R^n$, and so in the image, we are now in the position of Lemma~\ref{sublemma_2}, which says that, 
within the ball $U_v$, the portion of the plane $V_{\lambda} \cap U_v$ is mapped to $V_{\widetilde{\lambda}} \cap \varphi(U_v)$. Since this holds for every sufficiently small ball $U_v$, centered at any point $v \in V_{\lambda} \setminus W_{\lambda}$, and 
$V_{\lambda} \setminus W_{\lambda}$ is path-connected, we have $\varphi ( V_{\lambda} \setminus W_{\lambda}) \subset V_{\widetilde{\lambda}}$. Furthermore, as taking closures commutes under a homeomorphism, and $V_{\lambda} \setminus W_{\lambda}$ is dense in $V_{\lambda}$,
we in fact have that $\varphi(V_{\lambda}) \subseteq V_{\widetilde{\lambda}}$, which is what we needed to show. 
\end{proof}

Given a lattice point $\lambda \neq 0$, such that the open segment $(0, \lambda)$ contains no lattice points, define the collection $\mathcal{B}_{\lambda}$ of B-planes determined by $\lambda$ as follows
\begin{equation}
\mathcal{B}_{\lambda} = \{ V_{k \lambda} ~|~ k \in \Z \setminus \{0\} \}.
\end{equation}
The B-planes in the collection $\mathcal{B}_{\lambda}$ are parallel. Consequently, $V_{\lambda}$ separates $0$ from $V_{2 \lambda}$, the latter separates $V_{\lambda}$ from $V_{3 \lambda}$ etc. 
As the same separation properties hold for their images, since $\varphi$ is a homeomorphism, combining this with Lemma~\ref{planes}, it follows that $\varphi(\mathcal{B}_{\lambda}) = \mathcal{B}_{\widetilde{\lambda}}$ 
for a unique (modulo sign) $\widetilde{\lambda} \in \widetilde{\Lambda}$. We say that the collections $\mathcal{B}_{\lambda_1},...,\mathcal{B}_{\lambda_m}$ are linearly independent if $\nu(v)= m$, i.e. the dimension of the 
span of the vectors $\lambda_1,...,\lambda_m$ equals $m$. In what follows, we denote $\{ B_k \}_{k \in \N}$ and $\{ \widetilde{B}_k \}_{k \in \N}$ the Brillouin zones relative to $\Lambda$ and $\widetilde{\Lambda}$ respectively. 

\begin{lem}\label{lem_bound}
For every $k \geq 1$, we have that $\varphi (B_k) = \widetilde{B}_k$.
\end{lem}

\begin{proof}
Let $v \in \Int(B_k) \subset \R^n$ and let $w = \varphi(v)$. By Proposition~\ref{prop_focal_torus}, there are $k-1$ planes 
$ V_{\lambda_1},..., V_{\lambda_{k-1}}$ such that $V_{\lambda_s} \cap \ell_v \neq \emptyset$, $1 \leq s \leq k-1$ and $\mu(v) =0$. 
It follows that $\mu(w) = 0$. 

We need to show that $\iota(w) = k-1$. By Lemma~\ref{planes}, $\varphi (V_{\lambda}) = V_{\tilde{\lambda}}$, where $\lambda \in \Lambda$ and $\widetilde{\lambda} \in \widetilde{\Lambda}$.
As $\varphi$ is a homeomorphism for which $\varphi(0)= 0$, $V_{\lambda}$ separates $0$ and $v$ if and only if $V_{\tilde{\lambda}}$ separates $0$ and $w$ in $\R^n$. 
Hence $V_{\tilde{\lambda}_s} \cap \ell_w \neq \emptyset$, for every $1 \leq s \leq k-1$ where $V_{\tilde{\lambda}_s} := \varphi(V_{\lambda_s})$. Therefore, $\iota(w) = k-1$ and thus $w \in \widetilde{B}_k$. This proves that 
\[ \varphi( \Int (B_k)) = \Int(\widetilde{B}_k). \] 
By Proposition~\ref{prop_focal_torus}, we have $\cl( \Int (B_k) )= B_k$ and thus $\varphi (B_k) = \widetilde{B}_k$.
\end{proof}

\begin{lem}\label{lem_linearmap}
The homeomorphism $\varphi$ lies at a bounded distance from $A$, that is, 
\begin{equation}\label{eq_bound}
\varphi(v) = A(v) + \delta(v)
\end{equation}
with $A \in \GL(n, \R)$ and $\| \delta \|$ bounded.
\end{lem}

\begin{proof}
We will show that there exists a uniform tiling of $\R^n$ by identical parallelepipeds which is mapped by $\varphi$ to a tiling 
of $\R^n$ by the image parallelepipeds (which are again all identical). This will give us a linear map $A$ which also preserves this tiling. Because this tiling is uniform, 
the diameter, say $K$, of an image parallelepiped gives a uniform bound on the distance of a point $v \in \R^n $ under $\varphi$ and $A$, i.e.
\[ |\varphi(v)-A(v)| \leq K, ~\textup{for all}~v \in \R^n. \]
If we define $\delta(v):= \varphi(v)-A(v)$, then $|\delta(v)| \leq K$ and thus $\| \delta \|$ is bounded. 

Consider $\varphi (\mathcal{B}_{\lambda}) = \mathcal{B}_{\tilde{\lambda}}$ for an essentially unique $\tilde{\lambda} \in \widetilde{\Lambda}$, where $\pm \tilde{\lambda} \in \widetilde{\Lambda}$. 
As $\varphi(0) = 0$, we have that $\varphi( V_{k \lambda} ) = V_{\pm k \tilde{\lambda}}$. The collection $\mathcal{B}_{\lambda}$ does not consist of equally spaced B-planes, i.e. the distance 
between the two B-planes $V_{\lambda}$ and $V_{-\lambda}$ is twice the distance between any other two successive B-planes. However, 
\[ \mathcal{B}'_{\lambda} := \{V_{k\lambda} ~|~ k \in \Z\setminus \{0\}, ~k~\textup{odd}\} \subset \mathcal{B}_{\lambda} \] 
does consist of equally spaces B-planes. Choose ${\lambda}_1, {\lambda}_2, ...,{\lambda}_n \in \Lambda$ such that the corresponding collections $\mathcal{B}_{\lambda_i}$, with $1 \leq i \leq n$, are linearly independent. Then 
\begin{equation}\label{tiling_eq}
\R^n \setminus \bigcup_{s=1}^n \mathcal{B}'_{{\lambda}_s}
\end{equation}
defines a tiling of $\R^n$ by identical parallelepipeds. To finish the proof, we need to show that $\varphi$ maps this tiling to a uniform tiling of $\R^n$. It suffices to show that if we have $n$ linearly independent collections of B-planes as above, 
then the images of these collections under $\varphi$ are linearly independent as well. So let us prove this. Take $n$ linearly independent collections $\BB_{{\lambda}_1}, ..., \BB_{{\lambda}_n}$. The union of these B-planes tile $\R^n$ into $n$-dimensional parallelepipeds. 
Take one such parallelepiped $\PP$. Then $\partial \PP$ is topologically equivalent to $\mathbb{S}^{n-1}$. Because $\varphi$ is a homeomorphism, $\varphi (\partial \PP)$ is topologically equivalent to $\mathbb{S}^{n-1}$ as well. This can only be if the collections 
$\BB_{{\lambda}_1}',..., \BB_{{\lambda}_n}'$ are linearly independent. As the collections $\mathcal{B}_{\tilde{\lambda}_1}, ..., \mathcal{B}_{\tilde{\lambda}_n}$ are linearly independent, the union of the B-planes from the collections 
$\mathcal{B}'_{\tilde{\lambda}_1}, ..., \mathcal{B}'_{\tilde{\lambda}_n}$ form a tiling of $\R^n$. This tiling is preserved by $\varphi$ which now can be approximated by the linear map $A$ as described in~(\ref{eq_bound}).
\end{proof}

Let $\mathcal{A}(r,R) \subset \R^n$ be the Euclidean annulus with inner radius $r$ and outer radius $R$ centered at $0 \in \R^n$. The following result describes the asymptotic properties of the Brillouin zones, see~\cite{jones} and~\cite{skri}.
We give a sketch of the proof for the convenience of the reader.

\begin{lem}\label{lem_asymp}
For every $k$ there exist $0 < r_k < R_k < \infty$ such that 
\[ B_k \subset \mathcal{A}(r_k, R_k), ~\tu{with}~ | R_k - r_k | \leq C, \] 
for some constant $C>0$ independent of $k$.
\end{lem}

\begin{proof}[Sketch of the proof]
The condition that $v \in \Int(B_k)$ is, by Proposition~\ref{prop_focal_torus}, equivalent to the number of lattice points of $\Lambda$ contained in the Euclidean ball $D(v, |v|)$ being equal to $k$. As this number, up to an error relatively 
small compared to the total number of lattice points for large $|v|$, does not depend on $v$ but only on $|v|$, the claim follows for $\Int(B_k)$ and therefore for $B_k$ by passing to the closure.
\end{proof}

We now finish the proof. 

\begin{proof}[Proof of Theorem A]
By Lemma~\ref{lem_linearmap}, there exists a linear $A \in \GL(n, \R)$ such that $\varphi$ is homotopic to $A$, that is, there exists a uniformly bounded $\delta: \R^n \rightarrow \R^n$ such that $\varphi(v) = A(v) + \delta(v)$. 
By Lemma~\ref{lem_bound}, we have that $\varphi(B_k) = \widetilde{B}_k$ for every $k \geq 1$. By Lemma~\ref{lem_asymp}, the Brillouin zones $B_k$ and $\widetilde{B}_k$ are contained in Euclidean annuli with uniformly bounded thickness. 
It thus follows that, as $\delta$ is uniformly bounded, $A$ sends spheres to spheres. Therefore, by standard linear algebra, $A$ is an orthogonal mapping, i.e. $A \in \tu{O}(n,\R)$. To finish the proof, we must show that 
$\widetilde{\mathcal{V}} = \varphi( \mathcal{V} ) = A( \mathcal{V})$. From this it follows that $\widetilde{\Lambda} = A(\Lambda)$ as every $B$-plane is the perpendicular bisector of a point of the corresponding lattice.

To prove this, it suffices to show that $\varphi( \mathcal{B}_{\lambda} ) = A( \mathcal{B}_{\lambda} )$ for every $\mathcal{B}_{\lambda}$. We define the axis of $\mathcal{B}_{\lambda}$ to be the line through the origin perpendicular to the B-plane 
$V_{\lambda} \subset \mathcal{B}_{\lambda}$. Let $l$ be the axis $\mathcal{B}_{\lambda}$ and $\tilde{l}$ the axis of $\varphi(\mathcal{B}_{\lambda})$ and denote $l' = A(l)$. We claim that $\tilde{l} = l'$. Indeed, if $\tilde{l} \neq \l'$, then the planes 
$\varphi( V_{\lambda} )$ and $A( V_{\lambda})$ are not parallel, with $V_{\lambda} \subset \mathcal{B}_{\lambda}$. We can choose a point $v \in V_{\lambda}$ (with big distance from the axis $l$), such that the ball centered at $A(v)$ with radius $K$ misses the 
hyperplane $\varphi( V_{\lambda} )$ altogether, contradicting Lemma~\ref{lem_linearmap}. Therefore, $\tilde{l} = l'$. Similarly, the distance between two successive planes in $A(\mathcal{B}_{\lambda})$ and $\varphi(\mathcal{B}_{\lambda})$ 
has to be equal, as assuming otherwise is easily seen to yield again a contradiction with the uniform bound on $|\delta(v)| = | A(v) - \varphi(v) |$. Therefore, $\varphi( \mathcal{B}_{\lambda} ) = A( \mathcal{B}_{\lambda} )$ and since this holds for every 
$\mathcal{B}_{\lambda}$, it follows that $\widetilde{\mathcal{V}} = \varphi (\mathcal{V}) = A(\mathcal{V})$. 

Thus, as the lattices $\Lambda$ and $\widetilde{\Lambda}$ are related by an element of $\tu{O}(n,\R)$, the tori are isometric, up to rescaling. 
\end{proof}

\section{Concluding Remarks}\label{sec_questions}

The central question raised is what geometrical information of a closed Riemannian manifold is encoded in the mere topology of its focal decomposition.
We showed that for flat $n$-tori, the isometry type is essentially determined by the topology of the focal decomposition at a single tangent plane, since the torus admits a transitive group of isometries.

In general however, considering the focal decomposition at a single basepoint would not suffice in order to determine the geometry of the underlying manifold. Indeed, if one takes 
a canonical two-sphere and an ellipsoid of revolution, then the focal decomposition in two tangent planes, one in either surface, being homeomorphic, or even equal, does not imply these are isometric, 
as an ellipsoid of revolution has umbilical points at which basepoint the focal decomposition equals that of a canonical sphere. At other basepoints the focal decomposition of an ellipsoid will be different from that
of the canonical sphere. 

Let us further remark the following. From the focal decomposition at a given basepoint, one can determine how many geodesics of a given length start at the basepoint and end at another certain point. 
Besides this information, the focal decomposition also contains information about how these loci of points {\em interrelate} and this information may be important in the following respect. A notion similar to, but in a sense weaker than, the focal decomposition, 
is the {\em length spectrum of a manifold}, which records lengths of closed geodesics of a manifold with multiplicity. In the setting of flat tori, this is determined as follows. 
Let $\Lambda \subset \R^n$ be a lattice of rank $n$ and record the radii $\rho$ for which the (Euclidean) sphere $S_{\rho}$ of radius $\rho$ centered at the origin in $\R^n$ meets $\Lambda$ and how many lattice point meet this sphere. Even though 
this information determines for example the two-torus up to isometry, it fails to do so in higher dimensions. More precisely, E. Witt in~\cite{witt}, showed that, in dimension 16, there exist lattices $\Lambda, \widetilde{\Lambda} \subset \R^{16}$ that are not $\SO(16, \R)$ equivalent, yet have equal spectra. Using this information, J. Milnor in~\cite{milnor} remarked that the sequence of eigenvalues of the Laplace operator of a compact Riemannian manifold $(M,g)$ does not, in general, characterize its Riemannian metric. 

M. Peixoto in~\cite[page 461]{pei_3} discusses a problem similar to that of the length spectrum, as follows. Take again a flat torus and now record the radii, with multiplicity, for which the sphere $S_{\rho}$ in the cover either (1) is tangent to a B-plane or (2) meets an intersection of two or more B-planes. Call this sequence the {\em focal spectrum} of the torus. It is seen that the length spectrum above is a subsequence of this focal spectrum, but the focal spectrum records radii that are not present in the length spectrum. This focal spectrum is very sensitive to the position of the lattice points in space
and appears to be a finer measurement than the ordinary length spectrum. Let us pose here the question whether for flat $n$-tori, with $n \geq 2$, the focal spectrum determines the geometry of the torus up to isometry.

\subsection*{Acknowledgements}
The authors thank the referee for constructive and helpful comments and suggestions that helped improve the manuscript. Further, the authors thank Sebastian van Strien for several comments on the manuscript. 
Ferry Kwakkel received financial support from CNPq grant 301838/2003-8 and Marie Curie grant MRTN-CT-2006-035651 (CODY). Mauricio Peixoto received financial support from CNPq grant 301838/2003-8.

\end{document}